\documentclass{article}
\usepackage{amsmath,amssymb,amsfonts,amsthm,graphicx,enumerate}
\usepackage{graphicx}
\usepackage{color}
\usepackage{float}
\usepackage[english]{babel}
\usepackage{rotating}
\usepackage{caption}

\providecommand{\keywords}[1]{\textbf{\textit{Keywords:\ }} #1}

\theoremstyle{remark}

\newtheorem{Remark}{Remark}

\begin{document}

\title{Application of robust control for CSR formalization and stakeholders interest}

\maketitle

\author{Sana Ben Abdallah \footnote{Univ. Manouba, ESCT, Campus Universitaire Manouba, 2010, Tunisi; Address : Ecole Supérieure de Commerce de Tunis, La Manouba - Tunis - 2010 -Tunisie; E-mail: sana.benabdallah@esc.u-manouba.tn}
 \and Dhafer Saidane \footnote{Skema Business School, Université Cote d’Azur; Address: SKEMA Business School, Avenue Willy Brandt, 59777 Euralille –France; E-mail: dhafer.saidane@skema.edu} \and Mihaly Petreczky \footnote{Centre de Recherche en Informatique, Signal et Automatique de Lille (CRIStAL), UMR CNRS 9189, CNRS; Address: Ecole Centrale de Lille,59651 Villeneuve dAscq; E-mail: mihaly.petreczky@ec-lille.fr }}

\begin{abstract}
In this paper, we propose a new definition of sustainability that includes dynamics and equity. We propose a theoretical framework that allows finding a fair equilibrium and sustainable strategies for all stakeholders. The goal is to propose a strategy such that in the long run the attributes get close to an equilibrium point which is Pareto-optimal for the stakeholders. In order to calculate a strategy, we model stakeholders and actors as dynamical systems in state-space form. Furthermore, we use robust control and linear matrix inequalities (LMIs) to calculate the desired feedback strategy. We use several simulation scenarios to show the effectiveness of our proposed framework.
\end{abstract}
\keywords{ Sustainability, corporate social responsibility, stakeholder theory, control theory, linear matrix inequalities (LMIs), Pareto-optimality. }

  Declarations of interest: none

\newpage
\section{Introduction}

It is difficult for any company to satisfy the needs of all its stakeholders at the same time. In general, only the interests of the shareholders are protected by the managers. Those of other stakeholders are almost ignored  in particular because of their conflicting nature. The sustainable development and corporate social responsibility (CSR) approach tries to harmonize the opposing utilities of the different stakeholders. It defends a compromise between the expectations of these latter, whatever their power.

This paper proposes a reflection on the contemporary positioning of companies as "sustainable", and on a type of practices that have developed in recent decades, that of "dialogues" with the different "stakeholders".
CSR is defined by the European Commission (2002) as "a concept whereby companies integrate social and environmental concerns in their business operations and in their interaction with their stakeholders on a voluntary basis". It is therefore intended for all stakeholders. In general terms, organizations' engagement with stakeholders can be defined as the process of taking into account, the participation or involvement of individuals and groups who influence or are influenced by the company's activities.
Although opinions differ on how to integrate them, it seems that the integration of stakeholders is one of the essential ingredients for corporate sustainable excellence. All pose the question of the distribution of interests, resources, and responsibilities of each actor and highlight the strengths and limitations of stakeholder theory. It involves not only identifying, classifying and consulting stakeholders but also integrating their interests to advance the business and to be sustainable.

In this paper, we propose a new definition of sustainability that includes dynamics and equity. Indeed, we propose a theoretical framework that allows finding a fair equilibrium and sustainable strategies for all stakeholders.

Our approach is inspired by Rawls' theory. Indeed, this theory is founded on the idea that equity is deduced from the acceptability of the result. Traditional morality is replaced here by the agreement on mutually beneficial conventions. The rawlsian theory considers society as a system of cooperation accepted between free, equal and rational people. According to Rawls, the person, engaged in social cooperation is assumed to be capable of adjusting his objectives and aspirations according to what he can reasonably hope to obtain given his perspectives and position in society (\cite{rawls2005political}). The theory implicitly assumes that individuals are not slaves to their tastes and desires, but that they have some control over their preferences, and their effective choices: they are, in a certain sense, responsible for them. \cite{gauthier1986morals}  is the main philosopher who defends this thesis: he seeks to derive morality from rationality "non-coercive reconciliation of individual interest with mutual benefit". This is a form of justice as a mutual advantage, stemming from the Hobbesian tradition, where everyone seeks to best satisfy their interests, and, for this, accepts, rationally, constraints.

\paragraph*{The suggested approach}
We model the interests of stakeholders by utility functions which depend on attributes. These attributes represent various economic quantities, and they change in time. The change of attributes is modeled by a dynamical system in state-space form.  Space-state models (\cite{Moura}; \cite{Choi}; \cite{Souza}) have been widely used in economics. A state-space model describes how attributes are generated from their lag and how certain exogenous actions can influence those attributes. These exogenous actions represent the action of certain actors/stakeholders. The goal is to propose a strategy for choosing these exogenous actions, such that in the long run the attributes get close to an equilibrium point which is Pareto-optimal for the stakeholders: the equilibrium point cannot be changed without decreasing the utility function of one of the stakeholders. 
This means that in the long run, the strategy will lead to a \emph{fair} outcome,  as Pareto optimality  is used to express social justice between the stakeholders. Note that Pareto-optimality is closely related to Rawl's theory of social justice (\cite{gauthier1986morals}).  

Moreover, the strategy should be such that the vital interests of the stakeholders are not violated, i.e., at no time the attributes take a value for which the utility of one of the stakeholders descends below a certain level.
The situation where the utility function descends below a certain level represents a catastrophic scenario, whereby the corresponding stakeholder is forced to abandon the economic process. It could correspond to 
resource depletion, bankruptcy or basic needs not being made (if the stakeholder is a human being).
We call collections of attribute values \emph{sustainable}, if the utility functions of all stakeholders are above this critical threshold, when evaluated at these attribute values.

This property will be used to define sustainability. 
Sustainability means that we avoid situations where one of the stakeholder's vital interests are not respected. 
 By avoiding such situations we can ensure that
stakeholders will continue to cooperate and catastrophic outcomes will be avoided. 
This definition of sustainability captures the intuitive meaning of the concept: namely, that the economic process can be continued for a long period without a  major crisis.
That is, the strategy is both \emph{fair, acceptable and sustainable}.

Moreover, the proposed strategy is robustly sustainable i.e., even with the presence of disturbances or modeling errors, our strategy will still be sustainable.

Remarks are in order concerning the origins of the state-space models and the implementation of the strategy.

The state-space models can be viewed as arising from the behavior of rational agents acting using local information. 
More precisely, assume that each attribute corresponds to an agent who has the right to change them. Note that agents need not coincide with the stakeholders, they are two different concepts. 
Each agent  tries to choose the next attribute  according to that current attribute in such a manner that the utility function of the agent is maximized, assuming that all the other attributes remain constant. 
Note that the utility functions of agents are different from those of the stakeholders. 
It turns out that the equilibrium point which is Pareto-optimal from the point of view of stakeholders is a local Nash-equilibrium from the point of view of agents.

Concerning the implementation of the strategy, it can be done by imposing a suitable tax on the actors which are responsible for the exogenous actions. 

We use methods from robust control theory to calculate the described strategy.

Our contribution consists in the method used. To our knowledge this is the first study that attempts to find a fair and sustainable strategy that takes into account temporal aspects.


The structure of the remainder of the paper is as follows. The theoretical foundations' section will be devoted to the literature review; we will focus on the stakeholders' theory. Section 3 sets up the methodology. First, we will represent the general approach. In the next step, we will represent the calculation of the strategy and the safe invariant set. In section 4, a numerical case study will be performed to illustrate the proposed approach. First, we will present our phenomenological dynamic model. In the second subsection, we will summarize our results. Finally, an overall conclusion will be found in the last section.

\section{Related literature: Stakeholder theory}

The concept of sustainable development encourages companies to involve stakeholders in their governance. The issue of sustainable development is linked to the integration of the expectations and interests of stakeholders in corporate strategy and management \cite{Sharma}. As \cite{Capron} point out, "The concept of stakeholders is ubiquitous in all the literature on Corporate Social Responsibility".

The idea of only appealing to shareholders is now considered obsolete by several experts. The company, as part of a network of actors, must take into account the interests of its stakeholders. Indeed, according to \cite{Persais} : "One of the main challenges of the current leader is therefore to integrate the (non-economic) interests of a set of stakeholders and make them compatible with the interests of shareholders". Furthermore, \cite{Donaldson} argue that "all persons or groups with legitimate interests participating in an enterprise do so to obtain benefits and there is no prima facie priority of one set of interests and benefits over another", 

The stakeholder approach has been the subject of both empirical and theoretical studies. Today, stakeholders are at the heart of the social responsibility mechanisms implemented in companies. According to a broad consensus, stakeholder theory represents a relatively solid foundation, at least well established and recognized, for research on CSR, Business and society relations or business ethics. It is also used in debates on corporate governance and on the relationship between corporate strategy and sustainable development.

The concept of stakeholders is given by \cite{Freeman}. It is defined as "any group or individual who can affect or is affected by the achievement of the organization objectives." It is a concept that opens towards "a pluralist vision of the organization, an entity open to its environment” \cite{Samuel}. The stakeholder theory presents itself as an attempt to found a new theory of the firm integrating its environment to go beyond the traditional profit-making vision of the firm \cite{Gond} . Therefore, this theory seeks to integrate the interests of individuals and groups of people concerning the company and taking into account the social performance of this latter(\cite{Padioleau}).\cite{Sternberg} summarizes the concept of stakeholders by "Any person may have an interest in an organization". The definition of this concept is still the subject of many discussions (\cite{Pedersen}). Furthermore, the stakeholders' approach includes the views of stakeholders and makes them compatible with the views of shareholders. This is one of the most important challenges facing companies. 

Some authors have tried to classify the stakeholders in two visions. "Normative" or "Instrumental" (\cite{BALLET}). The normative vision is a purely ethical vision, where the company seeks to satisfy all stakeholders, by defining moral guidelines and use these guidelines as the basis for decision making. In contrast, the instrumental view is the consequence of taking into account stakeholders opinions as an essential element that leads to value creation. Managing relationships with stakeholders is a way for the company (directors and shareholders) to achieve its goals. \cite{Sharma} returned stakeholders into two groups: economic and non-economic. Economic stakeholders include all stakeholders involved in economic life and in productive activities of the company such as shareholders, suppliers, customers, etc. While the non-economic stakeholders associated with the environmental and the social actors. Also, they are linked to ethical dimensions. In short, stakeholders are defined as suppliers, customers, shareholders, employees, managers, regulators, and civil society ... etc (\cite{AVKIRAN}).

Some studies have also tried to prove the positive relationship between financial performance and the inclusion of stakeholder’s points of view (\cite{Luffman}; \cite{Jones}; \cite{HILLMAN}). Some other studies on stakeholder’s management also indicated a positive relationship between the plural form in management, the including of all stakeholders opinions, and the financial performance. For instance, \cite{Post} showed that among 89 studies, 48 of them showed this positive relationship. \cite{Tiras} argue also that companies that hold a good relationship with stakeholders exhibit higher performance. Thus, the integration of stakeholders can reduce risk; enhance the confidence of civil society, and improve the transparency of the regulatory framework (\cite{Holliday}). According to \cite{Sharma}, “In the short term, the integration of stakeholders can reduce costs and provide opportunities for differentiation. In the long term, it allows the dynamic construction of valuable competitive resources”. In most cases, we notice that effective stakeholder management enables banks to design policies for more efficient and stable banking systems (\cite{AVKIRAN2}).

 \paragraph*{The novelty of the paper concerning the related literature} 

The previous literature on stakeholder theory is focused only on the classification of the stakeholders and the impact of the inclusion of their points of view. To our knowledge no work has proposed a theoretical framework that allows guaranteeing the satisfaction of the entire stakeholders. 
In comparison with the previous discussed literature, we propose a model that considers the relationship of the companies with its stakeholders is equidistant and have the same importance. Our framework aims at transforming the conflictual system, formed of a group of individuals into an exchange relation with disparate objectives, to a group of individuals acting rationally in the name of a common objective. 

\section{Methodology}

As explained above, our goal is to find an acceptable strategy for all stakeholders i.e. with the choice of the actions, the state of the system will reach a certain sustainable state and remain in this sustainable state.

Mathematically, we represent the attributes of stakeholders as states of a dynamical system (\cite{DynSys}). A strategy is then a sequence of actions which drives the dynamical system to a certain state in which all the stakeholder are satisfied, i.e., their utility functions are above a certain threshold. Finding strategies for influencing dynamical systems is the core topic of control theory (\cite{Contr1}; \cite{franklin}). The proposed methodology consist of the following steps:


\begin{itemize}

\item Choice of the equilibrium point

In a first step we choose an equilibrium point in which the state of the system will reach a sustainable state, such that if the system is at that point, it will never leave it.
An equilibrium point, as known in the dynamical systems theory, is a state such as if the system reaches that state, then it will always remain there. The strategy we are looking for is one that forces the states of the dynamical system to approach the desired equilibrium point as time progresses. 

The state of the system may never become exactly the equilibrium point, it will get gradually closer and closer to it, hence the behavior of the system will get closer to its behavior in the equilibrium point. In particular, we will chose an equilibrium point where the utility functions of each stakeholder are above a certain threshold. We will call such equilibrium points sustainable. 

\item Calculating a safe set and a feedback strategy

In addition to reaching the equilibrium point it is necessary to find a set such that  it contains the equilibrium point and such that all the elements of this set are sustainable. Recall that by sustainability we mean that the utility function of each stakeholder is above a certain critical value. We will call such set, \emph{safe set}.
In parallel to calculating a safe set we also calculate a strategy such that when the strategy is applied, the safe set is \emph{invariant}. By invariance we mean that if the initial collection of attributes in this set, then at any time instance the attributes at that time instance will also be in that set. Moreover, under the application of this strategy, in the absence of disturbances, the attributes converge to the chosen equilibrium point. 

As a consequence,if we start in the safe set, we are sure that we will always remain there and converge towards the equilibrium point. If the initial state is not in the safe set, then the proposed strategy is not guaranteed to yield a sustainable behavior. Indeed, the safe set must ensure that all elements that belong to this set are sustainable and satisfy the constraints on attributes and actions. That is to say, even with the disturbance, the system will still be sustainable, although the attribute vector will no longer converge towards the equilibrium point.

\item Application of the strategy: feedback.
The calculated strategy will be in the form of feedback.  That is, at each time instance, the action prescribed by the strategy is a function of the current attribute values.The use of feedback and the properties of the safe set  guarantee that the strategy is robust. If the actual attribute values differ slightly from the ones prescribed by the model, due to external shocks (disturbances) or modeling error, but they are still in the safe set, then the application of the feedback will ensure sustainability and convergence to the equilibrium point in the absence of further disturbances. This property of feedback strategies is widely used in engineering (\cite{franklin}).

\end{itemize}




In the rest of this section, we will represent the general approach of our framework. Then, we will calculate the strategy and a safe invariant set. And we will end with a numerical case study in which, we will apply our phenomenological dynamic model.



\color{black}

\subsection{General approach}
\label{gen:approach}
 The idea is to model each attribute as a time varying variable, and 
 model the behavior of the stakeholders as a discrete-time state-space
 model (\cite{kailath}) of the form
 \begin{equation}
 \label{eq1}
   X(t+1)=F(X(t),U(t))
 \end{equation}
  where $X(t)=(X_1(t),\ldots,X_n(t))^T \in \mathcal{X} \subseteq \mathbb{R}^{n}$ 
  is the vector of attribute values at time $t=0,1,\ldots,$
  $U(t)=(U_1(t),\ldots,U_m(t))^T \in \mathcal{U} \subseteq \mathbb{R}^{m}$ 
  is the vector of actions by stakeholders,
  and $F: \mathcal{X} \times \mathcal{U} \rightarrow \mathcal{X}$ is the 
  \emph{state-transition} function. The function $F$ describes how the current attribute values 
  and the actions of the stakeholders influence the attribute values in the future. 
  The set $\mathcal{X}$ is the set of all possible values of the vectors of attributes and $\mathcal{U}$ is
  the set of all possible values of the vectors of actions by stakeholders.


For the purposes of this paper, we shall consider models where 
\begin{equation}
\label{state_space}
   \mathcal{X}=\{ (x_1,\ldots,x_n) \mid x_i \in [x_{i,min},x_{i,max}], i=1,\ldots,n\} ,
 \end{equation}
 and
\begin{equation}
\label{input_space} 
 \mathcal{U}=\{ (u_1,\ldots,u_m) \mid x_i \in [u_{i,min},u_{i,max}], i=1,\ldots,m\}. 
\end{equation}
That is,
the $i$th attribute is assumed to take values in the interval $[x_{i,min},x_{i,max}]$ and
the $j$th action is assumed to take values in $[u_{i,min},u_{i,max}]$. 

. 

Examples of actions $U(t)$ could be increase in minimal wage, or change in required solvency ratio, etc.

  Assume that we have a dynamical system of the form \eqref{eq1}. We call $X(t)$ the \emph{state} of
  \eqref{eq1} and we call $U(t)$ the \emph{input} of \eqref{eq1}. 
  Assume that there are $N$ stackholders and 
  for each stackholder there is an utility function   
   $f_i:\mathcal{X} \times \mathcal{U} \rightarrow [0,1]$, $i=1,\ldots,N$. 
  Intuitively, if $f_i(x,u)$ is close to zero, then the state and action pair $(x,u)$ is not favorable for the stackholder, if the value $f_i(x,u)$ is close to $1$, then
  the stakeholder is satisfied.


  We fix a set of values $\{f_{i,m}\}_{i=1}^{N}$ which represent the desired minima of the utility functions.

  In addition, we choose vectors $M  \in \mathbb{R}^{n+m}, m \in \mathbb{R}^{n+m}$ which denote the 
  maximum and minimum values of all attributes.

  We will call a state and input pair $(X(t),U(t)) \in \mathcal{X} \times \mathcal{U}$ \emph{sustainable}, if 
 \begin{equation}
 \label{eq:sus}
  \begin{split} 
    \forall i=1,2,\ldots,N:  f_i(X(t),U(t)) \ge f_{i,min},  \\
  \end{split}
  \end{equation}
  i.e., if in this state and input the value of the utility function of each stackholder is greater than 
  a certain minimal value. 
   In the sequel, we will concentrate on the case when 
  $U(t)$ is determined  a function of $X(t)$, i.e., $U(t)=\mathcal{F}(U(t))$ for some function
  $\mathcal{F}$. In this case, we say that the state $X(t)$ is sustainable, if
  $(X(t),U(t))$, $U(t)=\mathcal{F}(X(t))$  is sustainable. 

  Our goal is to find a strategy, i.e., a function $\mathcal{F}:\mathcal{X} \rightarrow \mathcal{U}$ such 
  that with the choice $U(t)=\mathcal{F}(X(t))$, the state of the system \eqref{eq1} will become sustainable.
  Moreover, we would like the strategy to yield be robustly sustainable, i.e.,
  in the presence of disturbances or modelling errors.

  To this end, we assume that $F(x,u)$ is affine, i.e. it is of the form
 \[ F(x,u)=\sigma_x(Ax+Bu+h)  \]
 where $A \in \mathbb{R}^{n \times n}$, 
$B \in \mathbb{R}^{n \times m}$ are suitable matrices, $h \in \mathbb{R}^n$ is a suitable vector, and $\sigma_x$ is a saturation function, i.e., 
\begin{equation}
\label{eq1.2satur} 
  \begin{split}
   & \sigma_x((x_1,\ldots,x_n)^T)=(\sigma_{1,x}(x_1),\ldots,\sigma_{n,x}(x_n))^T  \\
   & \sigma_{i,x}(x_i)=\left\{\begin{array}{rl} x_i & x_i \in [x_{i,min},x_{i,max}] \\
                                                x_{min} & x_i < x_{i,min} \\
                                                x_{max} & x_i > x_{i,max} 
                             \end{array}\right. 
\end{split}
\end{equation}
Note that if $z \in \mathcal{X}$, then  $\sigma_x(z)=z$, in particular, if $Ax+Bu+h \in \mathcal{X}$, then $F(x,u)=Ax+Bu+h$. 

 Hence, we assume that the dynamical system \eqref{eq1} takes the form
\begin{equation}
\label{eq1.2} 
  X(t+1)=\sigma_x(AX(t)+BU(t)+h).
\end{equation}

If we consider the equation \eqref{eq1.2} line by line, then the change in the value of the $i$th attribute is
\[ X_i(t+1)=\sigma_{i,x}(\sum_{i,j=1}^{n} a_{i,j} X_j(t)+\sum_{l=1}^{m} b_{i,l}U_l(t))
\]
\[ 
   A=\begin{bmatrix} 
      a_{11} & \ldots & a_{1n} \\
      a_{21} & \ldots & a_{2n} \\
      \vdots & \ldots & \vdots \\
      a_{n1} & \ldots & a_{nn} 
   \end{bmatrix}, ~  B=\begin{bmatrix} 
      b_{11} & \ldots & b_{1m} \\
      b_{21} & \ldots & b_{2m} \\
      \vdots & \ldots & \vdots \\
      b_{n1} & \ldots & b_{nm} 
   \end{bmatrix}
\]
Hence, if $a_{i,j}$ is positive (negative), it means that the increase in the value of the $j$h attribute leads to an increase (decrease)
in the value of the $i$th attribute in the next time step.
For example if $X_1(t)$ is profitability at time $t$, and $X_{10}(t)$ is the fixed wage at time $t$, then by increasing $X_{10}(t)$ we expect $X_1(t+1)$ to decrease (increase of wage leads to decrease of profitability), and hence $a_{1,10}$ should be negative

 In order to find a suitable strategy we will carry out the following steps. 
  \paragraph*{Choice of an equilibrium point}
     We find vectors $x_0 \in \mathcal{X}$, $u_0 \in \mathcal{U}$ such that 
     \begin{itemize}
     \item $x_0=Ax_0+Bu_0+h=F(x_0,u_0)$  (i.e. $(x_0,u_0)$ is an equilibrium point, that is
           if the $(x,u)$ is a solution of \eqref{eq1} such that $x(0)=x_0$ and
           $u(0)=u_0$ for all $t$, then $x(t)=x_0$.).
     \item $(x_0,u_0)$ is a sustainable state. 

     \item $f_i(x_0,u_0)=f_{i,t}$, $i=1,\ldots,N$ for some target values $f_{i,t} \ge f_{i,min}$ of the
      utility functions.  
     \end{itemize}

    That is, $(x_0,u_0)$ is such that if the system \eqref{eq1} is started in the initial state $x_0$
    and the input $U(t)$ is constant and it equals $u_0$, then the solution $X(t)$ will be equal to $x_0$.
    In other words, and equilibrium point is such that if the system is in that point, then it will never
    leave it.  Moreover, in the equilibrium point the utility functions take the target values $f_{i,t}$.

    The idea behind this is to choose $f_{i,t}$ in such a manner that all stakeholders are satisfied, e.g.,
    $f_{i,t}$ is larger than $0.5$ . In order to find $(x_0,u_0)$, 
    $x_0=(x_{0,1},\ldots,x_{0,n})^T$
    $u_0=(u_{0,1},\ldots,u_{0,m})^T$ the following non-linear programming problem should be
    solved:
    \begin{equation}
     \label{non_lin:prog:eq}
      \begin{split}
       & x_0=Ax_0+Bu_0+h \\
       & f_{i}(x_0,u_0) = f_{i,t}, ~ i=1,\ldots,N \\ 
       & x_{i,min} \le x_{0,i} \le x_{i,max}, ~ i=1,\ldots,n \\
       & u_{j,min} \le u_{0,j} \le u_{j,max}, ~ j=1,\ldots,m.
      \end{split} 
    \end{equation}

 \paragraph*{Choice of the strategy}
  Let us choose the strategy $\mathcal{F}$ as a feedback
  \begin{equation} \label{control_law} U(t)=\sigma_{u}(-K(X(t)-x_0)+u_0) \end{equation}
   where $K$ is a $m \times n$ matrix and
   \[
     \begin{split}
    & \sigma_u((u_1,\ldots,u_m)^T)=(\sigma_{1,u}(u_1),\ldots,\sigma_{m,u}(u_m))^T  \\
   & \sigma_{i,u}(u_i)=\left\{\begin{array}{rl} u_i & u_i \in [u_{i,min},u_{i,max}] \\
                                                u_{min} & u_i < u_{i,min} \\
                                                u_{max} & u_i > u_{i,max} 
                             \end{array}{rl}\right. 
    \end{split}
   \]
   That is, at every step, the input applied to the system depends on the current state.

  We would like to
   find $K$ and an ellipsoidal set $\mathcal{P}$ centered around $x_0$ of the form
    \begin{equation} \label{safety_region} 
    \mathcal{P}=\{ x \in \mathbb{R}^n \mid (x-x_0)^T Q^{-1} (x-x_0) < 1 \}
    \end{equation}
    where $P$ is an $n \times n$ strictly positive definite matrix, 
   such that the following conditions are satisfied:
   \begin{itemize}

   \item{\textbf{Stability}}
         If we use \eqref{control_law}, then $X(t)$ converges to $x_0$. Notice that if $X(t)$ converges to
         $x_0$, then $U(t)$ defined by \eqref{control_law} converges to $u_0$.

   \item{\textbf{Invariance}}
      If $X(0)$ belongs to $\mathcal{P}$ and $U(t)$ is chosen as in \eqref{control_law}, then
      $X(t)$ belongs to $\mathcal{P}$ for all $t$. 

   \item{\textbf{Safety}}
       If $X(t)$ belongs to $\mathcal{P}$ and $U(t)$ satisfies \eqref{control_law}, then
       for all $i=1,\ldots,N$, $(X(t),U(t))$ is a sustainable pair.

     \item{\textbf{Constraint satisfaction}}
     $\mathcal{P}$ should be a subset of $\mathcal{X}$ and for any $x \in \mathcal{P}$,
     $-K(x-x_0)+u_0 \in \mathcal{U}$.
   \end{itemize}

    That is, the matrix $K$ should be such that the application of the strategy \eqref{control_law}
    makes the system \eqref{eq1} \emph{stable} at the equilibrium point $(x_0,u_0)$, i.e. any solution
    $X(t)$  of \eqref{eq1} for the choice of $U(t)$ as in \eqref{control_law} is such that $X(t)$
    converges to $x_0$ and $U(t)$ converges to $u_0$. 
    Moreover, the set $\mathcal{P}$ is \emph{invariant} with respect to the system \eqref{eq1} and 
    the strategy \eqref{control_law}: if the strategy \eqref{control_law} is applied, and 
    $X(0)$ is in the set $\mathcal{P}$, then all the subsequent states $X(1),X(2), \ldots$ of
    \eqref{eq1} will be    in the set $\mathcal{P}$.  The set $\mathcal{P}$ is also \emph{safe}, i.e.,
    if a solution is in this set, then this solution is sustainable. 

    The motivation for constraint satisfaction is more involved. That $\mathcal{P}$ should be a subset
    of the set $\mathcal{X}$ of admissible states (attribute vectors) is not surprising, since all states 
    of interest live in $\mathcal{X}$. The reason that we require that $-K(x-x_0)+u_0 \in \mathcal{U}$
    whenever $x$ lies in $\mathcal{P}$ is the following: if this is not the case, then 
    $\mathcal{F}(x)=\sigma_{u}(-K(x-x_0)+u_0) \ne -K(x-x_0)+u_0$ and hence in this case were are not
    using the matrix $K$ for calculating our strategy. In turn, this may lead to instability and
    pathological behavior. This phenomenon is well known in control theory \cite{franklin}, and 
    our requirement for constraint satisfaction aims at avoiding this phenomenon. 

  \begin{Remark}[Robust sustainability]
    The strategy described above will be robust, if the initial state is in the safe set.
    More precisely,  if the true system is not $X(t+1)=F(X(t),U(t))$ but 
     $X(t+1)=F(X(t),U(t))+d(t))$  for some disturbance $d(t)$ such that $\|d(t)\|$ is sufficiently
     small, then with the strategy $U(t)=\sigma_u(-K((t)-x_0)+U_0)$ the system with perturbation $d(t)$
     will still be sustainable, although $X(t)$  will no longer converge to $x_0$. Indeed, 
     if $X(t) \in \mathcal{P}$, then $AX(t)+BU(t) \in \mathcal{P}$ and $U(t)=-K(X(t)-x_0)+u_0$ \footnote{If $x \in \mathcal{P}$, the $i$th component of $u=-K(x-x_0)+u_0$ belongs to $[u_{i,min},u_{i,max}]$}, 
     by the invariance property, and hence
     for small enough $d(t)$, $F(X(t),U(t))+d(t)$ will be in $\mathcal{P}$. 
That is, for small enough disturbances, if the state of the system is in $\mathcal{P}$, it will always remain there. Since the elements $x$ of $\mathcal{P}$ are sustainable (more precisely, $(x,u=-K(x-x_0)+u_0)$ is sustainable) it shows that 
     the proposed strategy is robustly sustainable.
  \end{Remark}

  \paragraph*{Summary}
    That is, we choose an equilibrium point $(x_0,u_0)$, and a strategy \eqref{control_law} and 
    a set $\mathcal{P}$ containing $x_0$, such that if  the initial state $X(0)$ belongs to the set $\mathcal{P}$, 
       it is also true that for all $t$, $(X(t),U(t))$ is sustainable.
      holds for all $t$. That is,if $\mathcal{P}$ is the set of sustainable initial
      states, such that if the system is started in such a sustainable initial state, then its state
     will always be sustainable.  This remains true even in the presence of small perturbation or modelling error.

    Moreover, in the absence of perturbation, attribute vector $X(t)$ will converge to $x_0$, so not only $(X(t),U(t))$ is sustainable, but eventually
    the value of the utility functions $f_i(X(t),U(t))$ will be close to $f_i(x_0,u_0)$. 

\paragraph*{Algorithm for calculating the desired strategy}
    In order to calculate the desired strategy, we propose to use tools from robust control, namely linear matrix inequalities (LMIs) \cite{LMIbook}. The details and the technical
    assumptions are described in Appendix \ref{lmi:appendix}.

\subsection{Pareto-optimality: Relationship with the classical approach}
\label{relationship_class}
\paragraph*{Interpretation of the result of feedback policy in terms of Pareto optimality}
  The proposed approach can be viewed as an attempt to achieve a Pareto-optimal (hence socially acceptable) outcome for  all stakeholders. More precisely, 
  we can choose the equlibirum point $(x_0,u_0)$ as follows:
  \begin{equation}
  \label{pareto:eq1}
   (x_0,u_0) = \mathrm{arg max}_{(x,u): F(x,u)=x, f_i(x,u) \ge f_{i,min}, i=1,\ldots,N} \sum_{i=1}^{N} f_i(x,u) 
  \end{equation} 
   This choice then guarantees that $(x_0,u_0)$ is a Pareto-optimal point for the utility functions $f_i(x,u)$, $i=1,\ldots,N$, i.e., it represents a socially desirable outcome.
   If the set of solutions of the constraints
   \begin{equation}
  \label{pareto:eq2}
        \mathcal{C}=\{ (x,u) \mid F(x,u)=x, f_i(x,u) \ge f_{i,min}, i=1,\ldots,N \} 
   \end{equation}   
 is not empty, then there is always a solution to \eqref{pareto:eq1}, if the utility functions $f_i$ are continuous, as $\mathcal{C}$ is a compact set.

\paragraph*{Implementation of the feedback policy via taxation}

 The strategy $U(t)$ then can be thought of as a policy to enforce the Pareto-optimal outcome. The strategy can be enforced by imposing a tax on the agents. Assume for the sake of simplicity that
 there is one agent, i.e., $m=1$, and the agent is the $N$th stakeholder. Note that the attributes change according to the dynamic equation $X(t+1)=F(X(t),U(t))$, so the various
 stakeholders are not able to change their attributes and hence the values of their utility functions. In this case, if we assume that the $N$th stakeholder chooses the value $U(t)$
 in such a manner that
 \begin{equation} 
 \label{pareto:eq3}
   U(t)=\mathrm{arg max}_{u_{1,min} \le u \le u_{2,min}} f_N(X(t),u) + I(X(t),u), 
 \end{equation} 
 where $I(X(t),u)$ is the tax to pay, then by  assuming that $f_N(X(t),u)$ is smooth n $u$ and $\dfrac{f_N(x,u)}{du}$ has constant sign for all $x \in \mathcal{X}$, $u \in [u_{1,min}, u_{1,max}]$, by choosing
 \[ I(x,u)=-\left(\frac{df_N}{du} (x,(-K(x-x0)+u_0))\right) u, \]
it follows that  $U(t)=\sigma_u(-K(X(t)-x_0)+u_0)$. That is, with a suitable choice of $I$, the optimal strategy of the agent will be to follow the control law $U(t)$. 

 That is, the calculated control law can be used by regulatory body to implement a policy via taxing, such that under this policy the attributes of the stakeholders converge to a Pareto-optimum.

\paragraph*{Justification of the state-space representation}
 Moreover, we assume that the function $F(x,u)$ describing the dynamics arises as follows:
 \begin{equation}
 \label{interp:model}
   F_i(x,u)=\sigma_{x,i}(x_i+\mu \dfrac{u_i(x,u)}{dx_i} )
 \end{equation}
where $u_i$ is an utility function which takes values in $[0,1]$, and
$\sigma_{x,i}$ is the saturation function from \eqref{eq1.2satur}.
 The interpretation of \eqref{interp:model} is as follows. The $i$th attribute depends on time, and \eqref{interp:model}
implies that
 \begin{equation}
 \label{interp:model2}
   X_i(t+1)=\sigma_{x,i}(X_i(t)+\mu \dfrac{u_i(X(t),U(t))}{dx_i})
 \end{equation}
 The interpretation of \eqref{interp:model2} is as follows. If we assume that each attribute belongs to an agent, then each agent will want to optimize its utility function $u_i$. More precisely,
 it wants to choose the next value $X_i(t+1)$ in such a manner that 
 $$u_i(X_1(t),\ldots,X_{i-1}(t),X_i(t+1),X_{i+1}(t),\ldots,X_n(t),U(t)) $$  
 is as large as possible. 
  That is, it assume that the
 attributes of the other agents will not change and it tries to adjust its own attribute value in such a manner that it optimizes its own utility function. However, 
 each agent has only limited information about its own utility function, and it can optimize it only locally, i.e., around $X_i(t)$. This means that the best agents can do is to
 change $X_i(t)$ in the direction which increases the utility the most. It is well-known that for some $\lambda > 0$, 
 \[
    \begin{split}
    & u_i(X(t),U(t)) \le \\
    & u_i(X_1(t),\ldots,X_{i-1}(t),
     X_i(t) +  \lambda \dfrac{u_i(X(t),U(t))}{dx_i},X_{i+1}(t),\ldots,X_n(t),U(t)), 
   \end{split}
 \] so $\dfrac{u_i(X(t),U(t))}{dx_i}$ is the direction into which 
 $X_i(t)$ should be changed in order to increase the utility function $u_i$. 

 Another way of looking at \eqref{interp:model2} is to assume that the time continuous values and  $X_i(\tau)$ satisfies $\dfrac{d}{d\tau} X_i(\tau) = \dfrac{u_i(X(\tau),U(\tau))}{dx_i}$.
 It then follows that 
 \[ \lim_{\tau \rightarrow +\infty} X_i(\tau)=\mathrm{arg max}_{x_i}  u_i(X_1(\tau),\ldots,X_{i-1}(\tau),x_i,X_{i+1}(\tau),\ldots,X_n(t),U(\tau)) \]
 under suitable technical assumptions. It then follows that
 $X_i(\tau+h)$ can be approximated by $X_i(\tau)+ \mu\dfrac{u_i(X(\tau),U(\tau))}{dx_i}$ for sufficiently small $\mu$ and if we identify $X_i(t)$ with $X_i(t\mu)$, \eqref{interp:model2} holds approximately, if we enforce the condition that $X((t+1)\mu)$ should belong to $\mathcal{X}$.

 Consider an equilibrium point $x_0=F(x_0,u_0)$ and assume that \eqref{interp:model} holds. Then $x_0=F(x_0,u_0)$ implies that
 $\dfrac{u_i(x_0,u_0)}{dx_i}=0$, i.e., the $i$th component $x_{0,i}$ of $x_0$ is the local optimum of $z \mapsto u_i(x_{0,1},\ldots,x_{0,i-1},z,x_{0,i+1},\ldots,x_{0,n},u_0)$. If 
 $\frac{d^2 u_i(x_0,u_0)}{dx_i^2} < 0$, then in fact $x_{0,i}$ is the local maximum o
 \[ z \mapsto u_i(x_{0,1},\ldots,x_{0,i-1},z,x_{0,i+1},\ldots,x_{0,n},u_0). \]

 This remark has the following implication: if we consider the game where each player corresponds to an attribute and each player tries to optimize the utility function $u_i$. The nodes
 of the game are attribute vectors $x \in \mathbb{R}^{n}$, the action of $i$th player is to choose the $i$th component of the next node. It then follows that if 
 $\frac{d^2 u_i(x_0,u_0)}{dx_i^2} < 0$ for all $i=1,\ldots,n$, then $(x_0,u_0)$ is a local Nash equilibrium. 

 Intuitively, this means that if we assume that each player can choose the $i$th component of the next node only locally, in the neighborhood of the current $i$th component, then the dynamics \eqref{eq1}
 arises as a repeated game, each node of the game corresponds to the state $X(t)$, and the next state $X(t+1)$ arises by each player trying to maximize its  utility function 
 $u_i$, i.e., $$X_i(t+1)=\mathrm{arg max}_{x_i \mbox{ close to } X_i(t)} u_i(X_1(t),\ldots,X_{i-1}(t),x_i,X_{i+1}(t),\ldots,X_n(t),U(t))$$. The equilibrium point $(x_0,u_0)$ is a Nash equilibrium of this game. 
 
 That is, the goal is to find a rule for choosing $U(t)$ such that the Nash equilibrium of this game such that this Nash equilibrium is also Pareto optimal w.r.t. utility functions $f_i$, $i=1,\ldots,N$. 
 That is, we have to types of utility functions:
 \begin{itemize}
 \item $u_i$, $i=1,\ldots,n$, are the utility functions which determine the interaction among various attributes and which define the mechanism of the time evolution of the attribute vectors,
 \item $f_i$, $i=1,\ldots,N$ are the utility functions of the stakeholders.
 \end{itemize}
 
 Note that the dynamics of the form \eqref{eq1.2} arise via the following choice of $u_i$:
 \begin{equation}
 \label{inter:model3}
 \begin{split}
  & u_i(x,u)=\frac{1}{\mu} \left( \frac{1}{2} (a_{ii}-1) x_i^2 + \sum_{j=1,i \ne j}^{n} a_{i,j} x_jx_i +  
 \sum_{j=1}^{m} b_{i,j}u_jx_i + h_i x_i\right)  \\
\end{split}
 \end{equation} 
 and $\frac{d^2 u_i}{dx_i^2}(x_0,u_0) < 0$, if and only if $a_{i,i} < 1$, i.e., if $u_i(x,u)$ is of the form \eqref{inter:model3} and $a_{ii} < 1$, then 
 $(x_0,u_0)$ is a local Nash equilibrium of the game above.

\section{Numerical case study}

\subsection{Phenomenological dynamic model}

In this example we will consider the behavior of stakeholders in the financial institutions (banks). However, our model can be applied not only to financing institution but in fact to any scenario of several rational stakeholders. The choice of this example was motivated by the previous researches in this field (\cite{SANA}).

In this case study, we use a phenomenological model to find the appropriate strategy. 
Our dynamic model is of the form \eqref{eq1}, more precisely, of the form \eqref{eq1.2}. 
That is, to define model, it is necessary to define  the state space $\mathcal{X}$, the input space $\mathcal{U}$,
the utility functions $f_i$ and the corresponding target values of  $f_{i,t}$ for $i=1,\ldots,N$.
and the equlibirum point $x_0,u_0$, and the matrices $A$ and $B$ and the vector $h$ from in \eqref{eq1.2}. 
We will only indicate the main steps, the precise numerical values will be presented in Appendix \ref{num_model}. 

\paragraph*{Choice of the state space and the input space }

In this example, we propose to the expert $6$ attributes for $3$ stakeholders.
The first stakeholder is the Manager. It is presented by two attributes: "Annual remuneration evolution" and "Return on assets (ROA)". Two attributes are also selected for the second stakeholder, the Regulator, that are, "Non-performing loans to total loans" and "Liquid assets to total assets". For the last stakeholder, the Customer, we also used two attributes: "Interest receivable to loans" and "Bank Fees to deposit".
The attribute "Bank Fees to deposit" will play the role of the input, as it is not influenced by the other attributes and can freely be set by the corresponding stakeholder ("Manager").

That is, for this example, 
$\mathcal{X} \subseteq \mathbb{R}^5$ is the state space and $\mathcal{U} \subseteq \mathbb{R}$ is the input space, more precisely,
$\mathcal{X}$ of the form \eqref{state_space} with $n=5$ and $\mathcal{U}$ is of the form \eqref{input_space} with $m=1$.
That is, each attribute and input are assumed to take values in the intervall $[x_{i,min},x_{i,max}], i=1,\ldots,5\}$ and $\mathcal{U}=[u_{1,min},u_{1,max}]$.
 All details are shown in Table \ref{table1} in Appendix \ref{num_model}.



\paragraph*{Determination of the utility function}
~~\\

The utility function $f_i(x,u$) is determined using multi-attribute utility approach due to \cite{Keeney}. This approach is based on the analyze of multiple variables simultaneously and assemble them on a synthetic indicator. 
In our case, the performance of each stakeholder is considered as a multi-attribute utility function (MAUF)  $(f_i(x,u))$ and each attribute is considered as a single attribute  utility functions (SAUF) $(g_{k}(x,u))$ i.e., we have 3 MAUF and 6 SAUF (2 SAUF for each stakeholder).

The first step consists in determining the attribute vectors ( i.e. $x_1, x_2, x_3, x_4, x_5$ and $u_6$). In general, these latter are assumed to be linear or exponential. According to \cite{Kim}, \textit {"when SAUFs are assumed to be linear or exponential, they are sufficient for most cases and their forms are solid"}. The SAUFs are determined using the ASSESS (to determine the three intermediate values) and LAB Fit (to fit the utility functions) software.

Equation \eqref{util_num} represents the forms of the single utility functions of the various attributes corresponding to risk-averse, risk-seeking and risk-neutral, respectively.

  \begin{equation}
  \label{util_num}
  \begin{split}
   & g_{k}(x,u)=\left\{\begin{array}{rl}
                   a-be^{(-cx)}  \\
                   a+be^{(cx)} \\ 
                   a +b(cx) \\
                   \end{array}\right. \\
  \end{split}
 \end{equation}
According to the expert's answers, the utility functions for each attribute are of the form:
 \begin{equation}
  \label{util_num2}
  \begin{split}
   & g_{k}(x,u)=\left\{\begin{array}{rl}
                   a_k+b_kc_k\frac{x_k}{100} &  ~ k \in \{1,2,4\} \\
                   a_k+b_ke^{c_k\frac{x_k}{100}} &  ~ k=3  \\
                   a_k+b_ke^{c_k\frac{u}{100}} &  ~ k=5 \\ 
                   a_k+b_ke^{c_k\frac{x_{k-1}}{100}} &  ~ k=6 
                   \end{array}\right. \\
 \end{split}
 \end{equation}
The scaling constants, the constants $a_k, b_k, c_k$, $k=1,\ldots,6$ are shown in Table \ref{table3}, Appendix \ref{num_model}.

To assess the utility function, a scaling constant $(k_{ij}, k_{i})$ is determined by the expert for each attribute and stakeholder, using ASSESS software, to establish the relevance of some with regard to other.

After determining the weights of the different attributes, we there by deduced the SAUF of the different attributes and consequently the MAUF of the different stakeholders.
The expert's opinion and our utility functions were obtained by using the same procedure and  data set used in \cite{SANA}. More details can be found in previous studies (see \cite{SANA} and \cite{Rebai})

Hence the MAUF ($f_i(x,u)$) of different stakeholders is of the following form:

\begin{equation}
  \label{util_num3}
  \begin{split}
   & f_i(x,u)=\frac{((K_ik_{2i-1}g_{2i-1}(x,u)+1)(K_ik_{2i}g_{2i}(x,u)+1)-1)}{K_i}, ~ i=1,2,3
  \end{split}
 \end{equation}
where $f_i(x,u)$ is the utility function for each stakeholder.$ k_i$ is the scaling constant for each attribute. $K$ is the overall scaling constant. 
The values of the constants $K_i,k_i$, $i=1,2,3$ are presented in Table \ref{table3}, Appendix \ref{num_model}.



 
%

\paragraph*{Choice of the matrices $A$ and $B$ of the model \eqref{eq1.2}}
We construct the matrices $A$ and $B$ of \eqref{eq1.2} as follows. First,
we propose to the expert a set of attributes to determine the effect of each attribute on the others. He will choose a value between $-1$ and $1$.\footnote{$-1$ if the attribute negatively affects the other attribute and $1$ if it positively affects it} The answers of the expert are then gathered in matrices $\hat{A}$ and $\hat{B}$.
Matrix $\hat{A}$ is composed of the attributes that can be influenced by the other attributes.
Matrix $\hat{B}$ is composed of attributes that are not influenced by any other attribute. We consider these attributes as inputs of the model.
The answers of the expert are presented in matrix $\hat{A}$ and $\hat{B}$  in Table \ref{model_table1} in Appendix \ref{num_model}.
From the discussion above it follows that the entries of $\hat{A}$ and $\hat{B}$ are numbers in the intervall $[-1,1]$ which express expert's opinion on the interaction between attributes. However, the values of the attributes belong to different intervals, hence the matrices $\hat{A}$ and $\hat{B}$ could be viewed as adequate models only for re-normalized attribute vectors, which take their values in the interval $[-1,1]$. Alternatively, the matrices $\hat{A}$ and $\hat{B}$ have to be rescaled in order to describe the dynamics of the true attribute vectors, hence the rescaling and shift by $b$. 
We therefore consider:  $A=0.5T^{-1}\hat{A}T$, $T^{-1}\hat{B}=B$ for a suitable diagonal matrix $T$.
The values of $A,B,T$ as indicated in Table \ref{mod1} in Appendix \ref{num_model}.

\paragraph*{Choice of the equilibrium point $(x_0,u_0)$ and the vector $h$ of \eqref{eq1.2}}

In order to find an equilibrium point $(x_0, u_0)$ we solve the following nonlinear programming problem:
 \begin{equation}
     \label{non_lin:prog:eq1}
      \begin{split}
       & f_{i}(x_0,u_0) = f_{i,t}, ~ i=1,\ldots,N \\ 
       & x_{i,min} \le x_{0,i} \le x_{i,max}, ~ i=1,\ldots,n \\
       & u_{j,min} \le u_{0,j} \le u_{j,max}, ~ j=1,\ldots,m.
      \end{split} 
    \end{equation}
using fin function of Matlab for  $f_{i,t}=0.7$, $i=1,2,3$.
Note that \eqref{non_lin:prog:eq1} differs from \eqref{non_lin:prog:eq}, as the first constraint of  \eqref{non_lin:prog:eq}, namely $x_0=Ax_0+Bu_0+h$ is absent from
\eqref{non_lin:prog:eq1}.  This is due to the difficulty assigning a value to $h$ in the absence of any data, as expert's opinion does not tell much about
the vector $h$. In fact, in some sense, the vector $h$ determines the equilibrium of the system: for any pair $(x_0,u_0)$, if we choose
$h=x_0-Ax_0-Bu_0$, then $(x_0,u_0)$ will be an equilibrium point for \eqref{eq1.2}.  Prompted by this observation, and by the lack of any other method to assign $h$
in the absence of measurement data, we propose to find first a candidate equilibrium point $(x_0,u_0)$ by solving \eqref{non_lin:prog:eq1}, and then
choose $h=x_0-Ax_0-Bu_0$. With this choice of $h$, $(x_0,u_0)$ will then be a true equilibrium point, in particular, it will be a solution of  \eqref{non_lin:prog:eq}.
While this approach is not satisfactory, we believe that it is justified for our purposes: our goal is to present a plausible example for illustration,
and we do not claim that our model is an adequate representation of some real economic process. For the real-life application of our approach, realistic models based
on econometric data will be required. In that case, the matrices $A,B$ and the vector $h$ will be estimated from real data.  
The values of $x_0,u_0$ and $h$ are indicated in Table \ref{mod1} in Appendix \ref{num_model}.



\paragraph*{Interpretation of the model using utility functions}
 As it was pointed out in Subsection \ref{relationship_class}, the dynamical system \eqref{eq1} can be interpreted as a result of rational behavior of 
 agents trying to optimize their own utility functions. More precisely, if we assume that there are as many agents as state components, the $i$th
 agent can influence the $i$th state component and at every time instance  it does so by choosing $X_i(t+1)$ in the vicinity of $X_i(t)$ so that
 $u_i(X_1(t),\ldots, X_{i-1}(t),X_{i}(t+1),X_{i+1}(t),\ldots,X_n(t),U(t))$ is maximal. Here $u_i$ is the utility function associated with the agent which manages
 the $i$th attribute. Then the right-hand side of \eqref{eq1} is of the form \eqref{interp:model}.
  If the dynamical system is of the form \eqref{eq1.2}, then this corresponds to the choice of $u_i$ of the form \eqref{inter:model3}.
 For the numerical example at hand the corresponding choices of $u_i$, $i=1,\ldots,5$ are presented in Table \ref{num:classical_inter}.

\subsection{Interpretation of the results}


In this section we present the results of our phenomenological dynamic model described in previous section. The strategy and simulations results are generated using Matlab.
We calculate the strategy using the method described in Appendix \ref{lmi:appendix}. The matrix $K$ defining the strategy of the form \eqref{control_law} and the matrix $Q$ defining the ellipsoid \eqref{safety_region}
can be found in Table \ref{num_result} in Appendix  \ref{num_model:calcul}. The details of the application of the method can be found in Appendix \ref{num_model:calcul}. 
As it was explained in Subsection \ref{relationship_class}, the strategy could be implemented using an appropriately formulated taxation policy.


 \paragraph*{Simulations}

 For illustration purposes, we show the simulation result of the system with the
  strategy \eqref{control_law}, if started from an initial state not in $\mathcal{P}$, and we also
  present the simulation result with the initial state being from $\mathcal{P}$, the values of the initial states are presented in
  Table \ref{initial_states:table}, Appendix \ref{num_model:simul}. 
  In both cases, the system we simulated was of the form
  \[
    \begin{split}
     & x(t+1)=\sigma_x(Ax(t)+Bu(t)+h)+d(t)) \\
     & u(t)=\sigma_u(-K(x-x_0)+u_0)
    \end{split}
  \]
  with
  \[ 
  \begin{split}
   & \sigma_x((x_1,\ldots,x_n)^T)=(\sigma_{1,x}(x_1),\ldots,\sigma_{n,x}(x_n))^T  \\
   & \sigma_{i,x}(x_i)=\left\{\begin{array}{rl} x_i & x_i \in [x_{i,min},x_{i,max}] \\
                                                x_{min} & x_i < x_{i,min} \\
                                               x_{max} & x_i > x_{i,max} 
                             \end{array}\right.  \\
   & \sigma_u((u_1,\ldots,u_m)^T)=(\sigma_{1,u}(u_1),\ldots,\sigma_{m,u}(u_m))^T  \\
   & \sigma_{i,u}(u_i)=\left\{\begin{array}{rl} u_i & u_i \in [u_{i,min},u_{i,max}] \\
                                                u_{min} & u_i < u_{i,min} \\
                                                u_{max} & u_i > u_{i,max} 
                             \end{array}{rl}\right. 
  \end{split}
\]
 The saturation functions $\sigma_x$ and $\sigma_u$ were used in order to make sure that the states and inputs stay in the sets $\mathcal{X}$ and $\mathcal{U}$ respectively.
 The latter was necessary because for states and inputs outside these sets, the utility function are not valid (their values no longer belong to the intervall $[0,1]$).

  The disturbance $d(t)$ represents the modelling error or external disturbances. While the calculation of the strategy was done
  for the case $d(t)=0$, since the strategy is stabilizing, it is robust with respect to small enough disturbances.  In fact, this is precisely the point of using a state-feedback. 
  We performed simulation with $d(t)=0$ and with the choice of $d(t)$ as in \eqref{choice:disturbance1}-\eqref{choice:disturbance2}, see Appendix \ref{num_model:simul} for a detailed explanation.

In order to better understand the results, we simulated the stakeholders' utility functions and the state components. Simulations are used to show the stakeholders' behavior in different situations (i.e., when the initial state in the safe set or not, when there are perturbations or not, etc.). In our case study we show the behavior of three stakeholders namely: the regulator, the client and the manager.

As presented in previous section, stakeholders are reasonably satisfied when the threshold of the utility function equals 0.7. From this threshold, we have deduced, using Matlab software, the equilibrium point of each attribute, which is defined in this case study as our initial state.

For Figures 1, 3, 5 and 6, the blue line represents the stakeholder's utility function and the red line represent the utility at equilibrium. For Figures 2 and 4, the blue line represents the state component and the red line represent the equilibrium value.

We start with the case where the system is started in the initial state is chosen from $\mathcal{P}$, and in the absence of perturbations (i.e. $d(t)=0$). Figures 1 and 2 below show the variation of the different stakeholders and attributes respectively.

\begin{figure}[H]
 \caption{Utility functions when the initial state is in $\mathcal{P}$ (see Table \ref{initial_states:table}), no perturbations ($d(t)=0$)}\label{fig:nodist:safe1}
\begin{center}
  \includegraphics[width=80mm,height=80mm,keepaspectratio, scale=0.8] {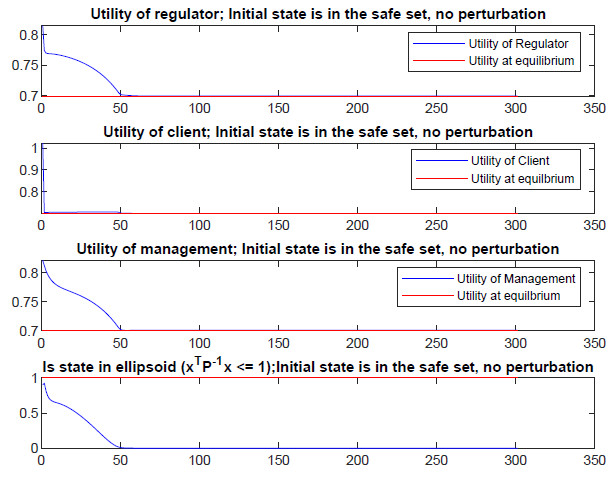}
\end{center}
  \end{figure}

\begin{figure}[H]
 \caption{State components and the input when the initial state is in $\mathcal{P}$ (see Table \ref{initial_states:table}), no perturbations ($d(t)=0$)  \label{fig:nodist:safe2}}
\begin{center}
  \includegraphics[width=80mm,height=80mm,keepaspectratio, scale=0.8]{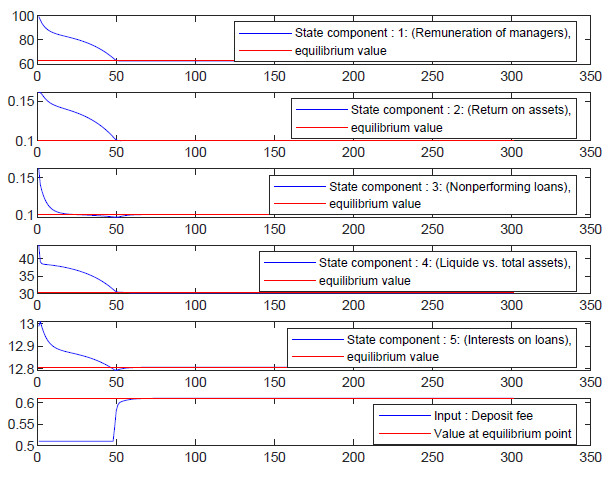}
\end{center}
  \end{figure}

As shown in Figures 1 and 2, we started with a sustainable utilities and states respectively. We see that the utility functions converge to the utility at equilibrium (i.e., 0.7) and the state components converge to the equilibrium values. Obviously, this situation guarantees us that the strategy will still be sustainable.

The second type of simulation can be seen from Figures 3 and 4. These figures show the utility functions and state components when we started 
with the initial state from $\mathcal{P}$ (see Table \ref{initial_states:table}), in the presence of perturbation ($d(t)$ as in \eqref{choice:disturbance2}). To deal with these disturbances, we have chosen a stabilizing feedback. 
\begin{figure} [H]
 \caption{Utility functions when the initial state is in  $\mathcal{P}$ (see Table \ref{initial_states:table}) and perturbations are present ($d(t)$ as in \eqref{choice:disturbance2})}\label{fig3dist2}
\begin{center}
  \includegraphics[width=80mm,height=80mm,keepaspectratio, scale=0.8]{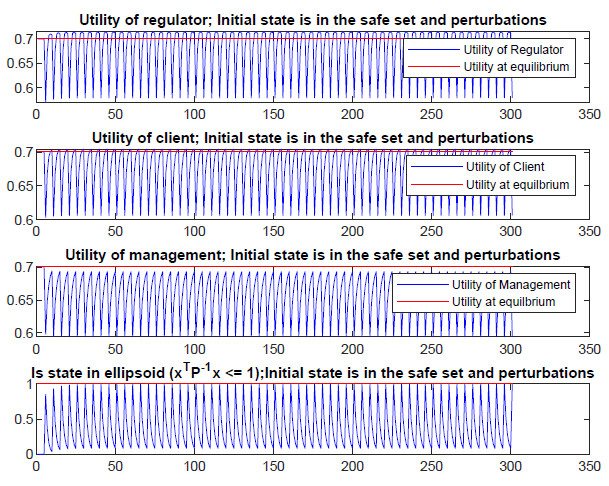}
\end{center}
\end{figure}

\begin{figure} [H]
\caption{State components and the input when the initial state is in  $\mathcal{P}$ (see Table \ref{initial_states:table}) and perturbations are present ($d(t)$ as in \eqref{choice:disturbance2}) }\label{fig3dist1}
\begin{center}
  \includegraphics[width=80mm,height=80mm,keepaspectratio, scale=0.7] {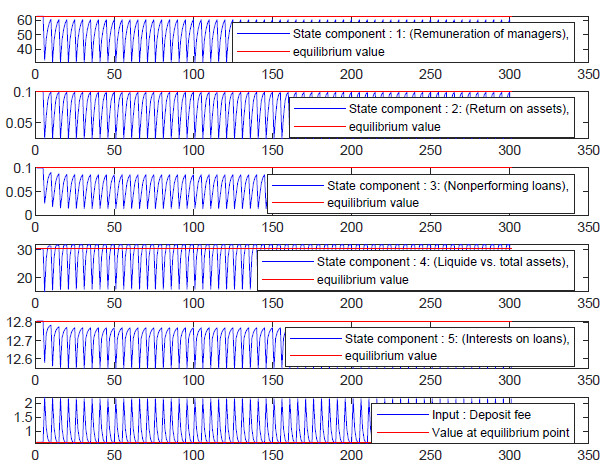}
\end{center}
  \end{figure}

Looking at Figures 3 and 4, we see a variation of the utility functions and the state components. We see that the proposed strategy is sustainable, and it remains so even in the presence of disturbance as long as the disturbance is small enough to keep the state in the safe set. In addition, we can conclude that, even with the disturbance, the system will still be sustainable, although the state components will no longer converge towards the equilibrium point.

Besides the simulations in the safe set (with and without disturbance), we also performed another simulation, but this time with an initial state that does not in the safe set and without disturbance (see Fig.5 below). In this case study, we obtained similar results as the previous ones, i.e., the proposed strategy reaches a sustainable level. But it is not guaranteed to remain in this sustainable behavior.

\begin{figure}[H]
  \caption{Utility functions when the initial state is not in $\mathcal{P}$ (see Table \ref{initial_states:table}), no perturbation ($d(t)=0$)}\label{fig1} 
\begin{center}
  \includegraphics[width=80mm,height=80mm,keepaspectratio, scale=0.8]{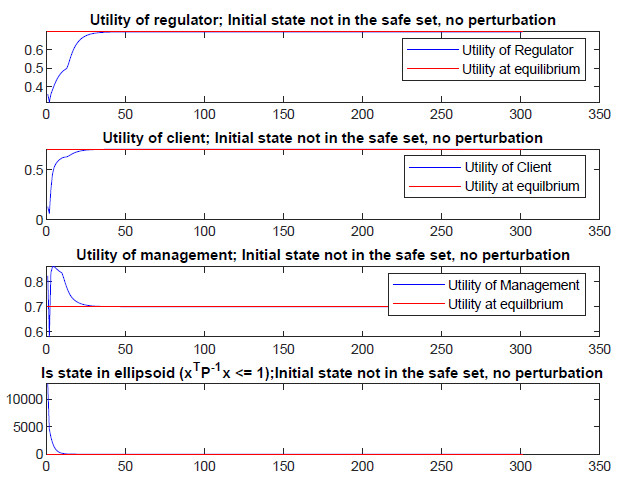}
\end{center}
  \end{figure}


Finally, we applied greedy input, i.e. we applied a strategy which is the best possible or the worst possible for one of the stakeholders. For our case study, we chose to increase the deposit fees for customers. That is, we chose the worst possible value for the customer.
Looking at Figure 6, we  notice that this action caused a remarkable decrease in the utility of the client (almost 0). (i.e that the client is not at all satisfied). In contrast, we note a remarkable increase of the utilities of the manager and regulator (almost 1) (i.e that the manager and regulator are totally satisfied). 
This action allows us to show the advantage of our model and the importance of the strategy found. A strategy that ensures the satisfaction of all stakeholders.

  \begin{figure}[H]
 \caption{Utility functions when the initial state is the equilibrium (see Table \ref{initial_states:table}), there is constant greedy input and no perturbation ($d(t)=0$)} \label{fig6}
\begin{center}
  \includegraphics[width=100mm,height=150mm,keepaspectratio, scale=0.7]{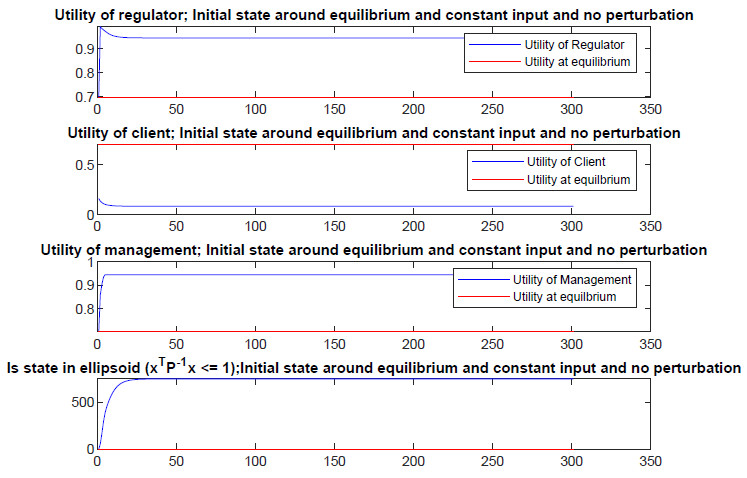}
\end{center}
  \end{figure}

\section {Conclusion}

The objective of this paper was to propose a theoretical framework that allows finding an acceptable strategy for all stakeholders by applying the control theory. To illustrate our approach, we used an academic example.
Note that our method requires a mathematical model in a state-space form that describes the interaction of various stakeholders. For the example at hand, we used a phenomenological model, which does not necessarily describe a real-life economic process.  The reason we did not use a more realistic model is that such a model should be  estimated from econometric data. The latter was not available for the example at hand. Even if the data had been available, building realistic models is a separate research topic, which would go beyond the intended scope of the paper.

Despite these shortcomings, we believe that the example of this paper demonstrated the feasibility of our method. 

Our approach helps us to correct some dangerous bad management practices, or deadly sins, that explain the current instability of the companies. Indeed, our model contribute to transforming the conflictual system, formed of a group of individuals with disparate objectives, to a group of individuals acting rationally in the name of a common objective. One of the main reasons for the unsustainability of companies is the excessive search for profit. Our approach consists of finding a strategy that allows everyone to reach an acceptable situation without looking for the ‘best’ situation. Furthermore, in general, stakeholders seek to maximize their profit against the other party. Our approach is to increase the claims of each stakeholder without exceeding a certain threshold. Finally, we propose a long-term vision instead of a short-term vision. Indeed, our strategy ensures that the system will reach a certain sustainable state and remain in this sustainable state.  

As previously stated, we applied the academic example used a phenomenological model. Future research should be directed to applying our methods to realistic models estimated from real-life data. 

\paragraph*{Practical implications}

The work that we are doing is the evaluation of strategies that allow the connection between stakeholders in the most beneficial way for all. This evaluation work can be performed by extra-financial rating agencies. Indeed, they can propose for companies a set of solutions (strategies) to ensure that all stakeholders are satisfied. This can improve their scoring process by giving additional notes, for example, on the application of good strategies. These scores could be viewed as a form (reputational) tax. Alternatively, tax policy could be adjusted to push stakeholders to apply good strategies. The construction of these strategies can then constitute the opportunity to stimulate a new form of negotiation that allows moving to a consensual and cooperative model, open to a truly sustainable economic environment.


\newpage


\appendix
\section{Calculating the strategy and a safe invariant set using linear matrix inequalities}
\label{lmi:appendix}
  In order to calculate the matrix $K$ and the set $\mathcal{P}$ described in Section \ref{gen:approach}, we
  use so called linear matrix inequalities (LMIs) \cite{LMIbook}. 
More precisely, we assume that the utility functions are piecewise quadratic and are of the form
\begin{equation}
\label{util3}
f_i(X,U)=\begin{bmatrix} X-x_0 \\ U-u_0 \end{bmatrix}_i^T Q_{i,k} \begin{bmatrix} X-x_0 \\ U-u_0 \end{bmatrix} + H_{i,k}\begin{bmatrix} X-x_0 \\ U-u_0 \end{bmatrix} + g_{i,k}, ~ \mbox{if}  (X-x_0,U-u_0) \in \mathcal{P}_{i,k}
\end{equation}
for suitably sized matrices $Q_{i,k}, H_{i,k}$ and scalar $g_i$, $i=1,\ldots,N$, $k=1,\ldots,D_i$, where
the sets $\mathcal{P}_{i,k}$ are polyhedral sets of the form
\begin{equation}
\label{util2.1}
 \mathcal{P}_{i,k}=\{ (x,u) \in \mathbb{R}^{n} \times \mathbb{R}^{m} \mid C_{i,k}\begin{bmatrix} x \\ u \end{bmatrix}+c_{i,k} \le 0 \}
\end{equation}
for suitable matrices $C_{i,k} \in \mathbb{R}^{r \times (n+m)}$ and vectors $c_{i,k} \in \mathbb{R}^{r}$, such that the union of all the sets
 $\mathcal{P}_{i,k} \cap \mathcal{X}$
covers the whole space $\mathcal{X}$.
\begin{Remark}
If the utility functions are not of the form \eqref{util3}, then they can be approximated with
arbitrary accuracy by piecewise-quadratic functions of the form \eqref{util3}. This follows from
the universal approximation property of piecewise-quadratic functions, which is a consequence of the universal approximation
property of piecewise-constant functions (which is a subclass of piecewise quadratic functions)/
\end{Remark}



We find  matrices $Q \in \mathbb{R}^{n \times n}$,  and $Y \in \mathbb{R}^{n \times m}$, $Q > 0$ \footnote{$Z > 0$, $Z < 0$ means that the matrix $Z$ is positive (negative) definite.}
,  by solving the following system of linear matrix inequalities (LMI)
 \begin{eqnarray}
  & \begin{bmatrix}
   - Q & QA^T-Y^TB \\ AQ-BY & -Q \end{bmatrix} < 0,  Q > W \label{lmi1} \\
 & \begin{bmatrix} Q & Qe_{i,n+m} \\ e_{i,n+m}^TQ & \frac{1}{\mu^2_i} \end{bmatrix} > 0, ~ i=1,\ldots,n  \label{lmi3.1} \\
  & \begin{bmatrix} Q & e_{j}Y \\  (e_jY)^T & \frac{1}{\mu_{n+j}^2} \end{bmatrix} > 0, ~ j=1,\ldots,m \label{lmi3.2} 
 \end{eqnarray}
  where $W \in \mathbb{R}^{n \times n}$ and $\mu_i$, $i=1,\ldots,n+m$ are design parameters chosen by the user, and
  $x_{0,i}, u_{0,j}$, $i=1,\ldots,n$, $j=1,\ldots,m$ denote the $i$th and $j$th component of $x_0$ and $u_0$ respectively.
  The parameters   $\mu_i$, $i=1,\ldots,n+m$ are chosen so that 
  if  $(x_i-x_{0,i} )^2 \le \mu_i^2$, then $x_i \in [x_{min,i},x_{max,i}]$ for all $i=1,\ldots,n$, and  if 
  $(u_j-u_{0,j})^2 \le \mu_{n+j}^2$, then $u_j \in [u_{min,j},u_{max,j}]$, $j=1,\ldots,m$. 

  This can be achieved by choosing $\mu_i \le \min\{x_{max,i}-x_{0,i}, x_{0,i}-x_{min,i}\}$, $i=1,\ldots,n$, and 
  $\mu_{n+j} \le \min\{u_{max,j}-u_{0,j}, u_{0,j}-u_{min,j}\}$, $j=1,\ldots,m$.

  The matrix $W$ should be chosen as a symmetric positive semi-definite matrix and it can be used to control the size of the ellipsoid $\mathcal{P}$. 

  We then set $K=YQ^{-1}$ and 
  \begin{equation}
  \label{ellipsoid2}
  \mathcal{P}=\{x \in \mathbb{R}^n \mid (x-x_0)^T Q^{-1} (x-x_0) < 1\}. 
  \end{equation}

  Finally, we verify that the following LMI with the indeterminate 
  $\tau > 0$, $\tau_{i,k,l} > 0$, $i=1,\ldots, N$, $k=1,\ldots,D_i$, $l=1,\ldots,r$ has a solution
 \begin{equation}
 \begin{split}
& \begin{bmatrix}
     \mathcal{S}^TQ_{i,k}\mathcal{S}  &  0.5 (H_{i,k}\mathcal{S})^T \\
     0.5 H_{i,k}\mathcal{S} &  g_{i,k}-f_{i,min} 
    \end{bmatrix} +  \tau \begin{bmatrix} Q^{-1} & 0 \\ 0 & -1 \end{bmatrix} + \\
    & + \sum_{l=1}^{r} \tau_{i,k,l} \begin{bmatrix} 
     0 & (e_{l,r}^TC_{i,k}\mathcal{S})^T0.5 \\
     (e_{l,r}^TC_{i,k}\mathcal{S})0.5 & e_l^Tc_{i,k} 
     \end{bmatrix} >0, 
  \label{lmi2} 
 \end{split}
 \end{equation}
 where $\mathcal{S}=\begin{bmatrix} I_n  \\ -K \end{bmatrix}$, and 
 $e_{i,d}$ is the $i$th standard unit vector of $\mathbb{R}^{d}$, i.e., all the elements of $e_i$ are zeros, except the $i$th one, which is $1$. 

The intuition behind the equations \eqref{lmi1} --  \eqref{lmi2} is the following.

  \begin{enumerate}
     \item  LMI \eqref{lmi1} ensures that the feedback $U(t)=-K(X(t)-x_0)+u_0$ will stabilize the system
    $X(t+1)=AX(t)+BU(t)+h$,    $\lim_{t \rightarrow \infty} X(t)=x_0$.

    \item LMI \eqref{lmi3.1} ensures that if $x \in \mathcal{P}$, then
         the $i$th component $x_i$ of $x$ satisfies $x_i \in [x_{i,min},x_{i,max}]$.
         Likewise, \eqref{lmi3.2}  ensures that if $x \in \mathcal{P}$ an $u=-K(x-x_0)+u_0$
         then the $j$th component $u_j$ of $u$ satisfies $u_j \in [u_{j,min},u_{j,max}]$.

%

   \item LMI \eqref{lmi2} ensures that
         if $x \in \mathcal{P}$ and $u=-K(x-x_0)+u_0$ then $f_{i}(x,u) \ge f_{i,min}$. 

       
  \end{enumerate}
   To sum up, the matrix $Q$ is calculated so that all the elements
   of the ellipsoid $\mathcal{P}$ are sustainable and satisfy the constraints on the
    attributes and actions.  

  From classical results of control theory it then follows that the matrix $K$ and the set $\mathcal{P}$ 
  satisfies the conditions described in the previous section. 
The solution of \eqref{lmi1}-\eqref{lmi2} is calculated using classical numerical tools YALMIP and its interface with Matlab. 


\section{Numerical example}
\label{num_model}

\subsection{Tables with the parameters of the example}
\begin{table}[H]
\caption{Selected Attributes \label{table1}}
\begin{tabular}{|l|l|l|}
\hline
Stackholders & Attributes & Range 
\\ \hline 
Managers:  
& Attribute $x_1$ : Annual remuneration evolution ($x_1$) & [-0.67 , 1.3] \\ 
& Attribute $x_2$: Return on assets (ROA) ($x_2$) & [-0.02 , 0.01] \\ \hline 
Regulator:  
& Attribute $x_3$: Non-performing loans to total loans ($x_3$) & [0.15 , 0]  \\
& Attribute $x_4$: Liquid assets to total assets ($x_4$) & [0.1 , 0.64]  \\ \hline 
Customers 
& Attribute $x_5$: Interest receivable to loans ($x_5$) & [0.31 , 0.015] \\ 
& Attribute $u$: Bank Fees to deposit ($u_1$) & [0.14 , 0.005] \\ \hline

\end{tabular}
\end{table}

\begin{table}[H]
\caption{Matrices $\hat{A}$ and $\hat{B}$ originating from expert's opinion \label{model_table1}}
\begin{tabular}{|c|c|}
\hline
  $\hat{A}$ & $\hat{B}$ \\
\hline
   $\begin{bmatrix} 1  &  0.8  & -0.2  &  0.5  &  0.2 \\
    0.2   & 1  & -0.4  &  0.6 &  0.5 \\
         0  & -0.4 &   1    &     0   & 0.7 \\
    0.2  &  0.5 &  -0.8   & 1  &  0.4 \\
         0  &  0.3 &   0.5 &   0.5   & 1
   \end{bmatrix}$ & 
$\begin{bmatrix} 
0.3 \\ 0.2  \\ 0\\ 0.2 \\ 0.4 \end{bmatrix}$ \\
\hline
 \end{tabular}
\end{table}



\begin{table}[H]
\caption{State-space transformation $T$ and the model parameters $A=0.5T^{-1}\hat{A}T$, $T^{-1}\hat{B}=B$, $h$, equilibrium $x_0,u_0$ \label{mod1}}
\begin{tabular}{|c|c|}
\hline
$T$ & 
 \( \begin{bmatrix}  0.0248      &   0     &    0     &    0   &      0    &     0 \\
         0  &     10     &    0    &     0    &     0    &     0 \\
         0  &       0 &      10   &      0    &     0     &    0 \\
         0   &      0  &       0   & 0.049  &      0       &  0 \\
         0   &      0   &      0     &    0  &  3.3333     \\    
         0   &      0  &       0    &     0    &     0   & 0.9487
     \end{bmatrix} \) \\
\hline
  $A$ & \( \begin{bmatrix}
   0.5  &161.4498  &-40.3624  &  0.4949  & 13.4541 \\
    0.0002 &   0.5   &-0.2  &  0.0015  &  0.0833 \\
         0  & -0.2 &   0.5   &      0 &   0.1167 \\
    0.0505  & 50.9759 & -81.5615  &  0.5  & 13.5936 \\
         0  &  0.45  &  0.75  &  0.0037  &  0.5
\end{bmatrix}  \) \\
\hline
   $B$ & \( \begin{bmatrix} 
12.1087 \\   0.02  \\     0 \\   4.0781  \\  0.12
 \end{bmatrix} \)
\\ 
\hline
$x_0$ & \(  \begin{bmatrix}  62.64 & 0.1 & 0.1 & 30.25 & 12.8            
\end{bmatrix}^T \) \\
\hline
$u_0$ &
      \( 0.61 \)
\\
\hline 
$h$ & $\begin{bmatrix} -175.43 & -1.07 & -1.42 & -161.54 & 6.1 \end{bmatrix}^T$ \\
\hline
\end{tabular}
\end{table}



  

\begin{table}[H]
\footnotesize
\centering
\caption{Utility functions for agents managing each attribute  \label{num:classical_inter}}
  \begin{tabular}{|c|c|}
\hline
    $u_1(x,u)$ & $-0.5 \cdot 0.5x_1^2+161.4498 x_2x_1-40.3624x_3x_1 + 0.4949x_4x_1+13.4541x_5x_1+12.1087u x_1 -175.43x_1$ \\ \hline  
    $u_2(x,u)$ & $-0.5 \cdot 0.5x_2^2+  0.0002x_1x_2-0.2x_3x_2 + + 0.0015 x_4x_2  +  0.0833 x_5x_2+0.0200u x_2-1.07x_2$ \\ \hline 
    $u_3(x,u)$ & $-0.5 \cdot 0.5 x_3^2-0.2x_2x_3+0.1167x_5x_3-1.42x_3$  \\ \hline 
    $u_4(x,u)$ & $-0.5 \cdot 0.5 x_4^2+0.0505 x_1x_4 + 50.9759 x_2x_4 -81.5615 x_3x_4  + 13.5936 x_5x_4+4.0781u x_4-161.54x_4$ \\ \hline
    $u_5(x,u)$ & $-0.5 \cdot 0.5 x_5^2+0.45x_2x_5 + 0.75x_3x_5 +  0.0037x_4x_5+0.12ux_5+6.1x_5$ \\ \hline
 \end{tabular}
\end{table}
\color{black}


\begin{sidewaystable}
\footnotesize
\centering
\caption{Scaling constant and and constants \label{table3}}

\begin{tabular}{| l | p{5cm} | l | l | l | l | l |}
\hline
Stakeholders  &  Attributes & $ K_i$  & $ k_{ij}$ &  $a_k$  & $b_k$  &  $c_k$
\\ \hline 
Managers  
& Attribute $x_1$: Annual remuneration evolution & -0.35854   &   0.59   &  0.3535 &  0.7625  & 0.7625  \\ 
& Attribute $x_2$: Return on assets (ROA)            &                       & 0.52     & 0.597   & 5.686      & 5.686  \\ \hline 
Regulator  
& Attribute $x_3$: Non-performing loans to total loans & -0,409837  &  0.52 & -0.2919 &  1.276 & -9.923 \\
& Attribute $x_4$: Liquid assets to total assets &    &    0.61  & 0.1441  & 1.342  & 1.342 \\ \hline 
Customers 
& Attribute $x_5$: Interest receivable to loans & -0.333333  &   0.5   &   -0.08086    &  1.172  &  -8.528      \\ 
& Attribute $u$: Bank Fees to deposit          &  &  0.6     &   -0.428   & 1.461      &  -8.653    \\ \hline

\end{tabular}
\end{sidewaystable}



\begin{sidewaystable}
\footnotesize
\centering
\caption{Constants  \label{table3.1}}



\begin{tabular}{|  l | p{3cm} | l | l | l | l | l | }

\hline
Stakeholders  &  Attributes & $x_{k,min}$  &   $s_{k,2}$ &  $s_{k,3}$  &  $s_{k,4}$  &  $x_{k,max}$
\\ \hline 
Managers  
& Attribute $x_1$: Annual remuneration evolution &   -67 \%   &  10\%   &  48\% & 79\%  & 103\%  \\ 
& Attribute $x_2$: Return on assets (ROA)  &     -2 \%   &  -0.24\%     & 0.16\%   &  0.64\%  & 1\%  \\ \hline 
Regulator  
& Attribute $x_3$: Non-performing loans to total loans &   15\%  & 8\% & 5.5\% &  1.7\% &  0\% \\
& Attribute $x_4$: Liquid assets to total assets   & 9.86\%   &   20.5\%  & 30.5\%  & 50\%  &  63.91\%  \\ \hline 
Customers 
& Attribute $x_55$: Interest receivable to loans &        31.05\%  &  14.5\%    &   9\%    &  3\%  &  1.52\%      \\ 
& Attribute $u$: Bank Fees to deposit     &     14\%  &  9\%     &   5.5\%   &  2\%     &  0.51\%    \\ \hline

\end{tabular}
\end{sidewaystable}


\subsection{Calculating a strategy}
\label{num_model:calcul}

 We would like to apply the method of Appendix \ref{lmi:appendix} to the numerical example. However, to this end, we have to solve a small technical issue.
Namely, not all the utility functions are of the form \eqref{util3}.  More precisely, for $i=1$, the utility function $f_{1}(x,u)$ is of the form \eqref{util3}, with
$D_1=1$, $\mathcal{P}_{1,1}=\mathbb{R}^{n+m}$ ($C_{1,1}=0,c_{1,1}=0$), and 
with the following choice of $Q_{1,1},H_{1,1},g_{1,1}$:
 \[ 
   \begin{split}
    Q_{1,1}=10^{-03} \begin{bmatrix} 
         0 &  -0.1034   &      0    &     0  &       0   &      0 \\
   -0.1034 &        0  &       0    &     0   &      0    &     0 \\
         0 &        0  &       0  &       0   &      0  &       0 \\
         0 &        0   &      0  &       0    &     0    &     0 \\
         0     &    0     &    0     &    0       &  0      &   0 \\
         0     &    0     &    0     &    0      &   0     &    0
   \end{bmatrix}, \\
  H_{1,1}=\begin{bmatrix} 0.0030   & 0.1555      &   0      &   0   &      0   &      0 \end{bmatrix} \\
  g_{1,1}=0.4958  
  \end{split} 
\]
However, for $i=1,2$, the utility functions $f_i$ are not of the form \eqref{util3}. In order to be able to apply Appendix \ref{lmi:appendix}, 
for $i=2,3$, the utility functions $f_i$ will be approximated by functions $f_{i,pl}$, $i=2,3$ of the 
form  \eqref{util3}. 
This will be done as follows. 
 For $i=2,3$, the utility functions $f_i$ can easily be approximated by functions of the 
 form  \eqref{util3} as follows. 

For $k=3,5,6$, let  $x_{k,max}=s_{k,1} > s_{k,2} > s_{k,3} > s_{k,4} > s_{k,5}=x_{k,min}$ be such that $g_k(s_{k,j})=0.25j$, $j=1,2,3,4,5$. 
For each $k=3,4,5,6$, the utility functions $g_{k}$ are approximated by piecewise-linear functions
$g_{k,pl}$
 \begin{equation}
\label{pl:approx:eq-1}
   \begin{split}
  &   g_{k,pl}(x,u)=n_{k,j}^T \begin{bmatrix} x \\ u \end{bmatrix} + b_{k,j}
     ~\mbox{ if }  R_{k,j}\begin{bmatrix} x \\ u \end{bmatrix} + r_{k,j} \le 0, ~ j=1,\ldots,4 \\
  & n_{k,j}=\frac{0.25}{s_{k,j+1}-s_{k,j}} E_k \\ 
  & b_{k,j}=\frac{0.25s_{k,j}}{s_{k,j+1}-s_{k,j}} \\ 
  & E_{4}=e_{4,n+1}^T, E_3=e_{3,n+1}^T, E_5=e_{6,n+1}^T, E_6=e_{5,n+1}^T \\
  & R_{k,j}=\begin{bmatrix} E_k \\ -E_k \end{bmatrix}, r_{k,j}=\begin{bmatrix} s_{k,j} \\ s_{k,j+1} \end{bmatrix}
 \end{split}
 \end{equation}
  The vectors $n_{k,j},b_{k,j},R_{k,j},r_{k,j}$ can readily be computed using the values $s_{k,j},s_{k,j+1}$
  in Table \ref{table3.1}, $k=3,4,5,6$, $j=1,\ldots,4$. 
 We then approximate the functions $f_2,f_3$ by the following functions
 \begin{equation}
 \label{pl:approx}
    \begin{split}
     f_{i,pl}(x,u)=\frac{((K_ik_{2i-1}g_{2i-1,pl}(x,u)+1)(K_ik_{2i}g_{2i,pl}(x,u)+1)-1)}{K_i}, ~ i=2,3
   \end{split}
 \end{equation}
 where $K_i,k_i$, $i=2,3$ are the same as in \eqref{util_num3}. 
 It then follows that $f_{2,pl}$ satisfies \eqref{util3} with $D_i=4$ and for all $j=1,\ldots,4$,
 \begin{equation}
 \label{pl:approx:eq1}
    \begin{split}
     & Q_{2,j}=\frac{N_4N_{3,j}^T}{K_2} \\
     & H_{2,j}= \frac{d_{3,j}N_4+d_4N_{3,j}}{K_2}+2\begin{bmatrix} x_0 \\ u_0 \end{bmatrix}^T Q_{2,j} \\
     & g_{2,j}=\frac{d_{3,j}d_4-1}{K_2}-\begin{bmatrix} x_0 \\ u_0 \end{bmatrix}^T Q_{2,j}\begin{bmatrix} x_0 \\ u_0 \end{bmatrix} + H_{2,j}\begin{bmatrix} x_0 \\ u_0 \end{bmatrix}\\
    & N_{3,j}=K_2k_3n_{3,j}^T, ~d_{3,j}=K_2k_3b_{3,j}+1 \\
    & N_4=b_4c_4K_2k_4, ~ d_4=a_4K_2k_4+1 \\
    & C_{2,j}=R_{3,j}, ~ c_{2,j}=r_{3,j}+R_{3,j}\begin{bmatrix} x_0 \\ u_0 \end{bmatrix}.
   \end{split}
 \end{equation}
Similarly $f_{3,pl}$ satisfies \eqref{util3} with $D_i=16$ and for all $j=1,\ldots, 16$, 
$j=4(j_1-1)+j_2$, $j_1,j_2=1,\ldots,4$,
\begin{equation}
 \label{pl:approx:eq2}
    \begin{split}
     & Q_{3,j}=\frac{N_{5,j_1}N_{6,j_2}^T}{K_3} \\
     & H_{2,j}= \frac{d_{5,j_1}N_{6,j_1}+d_{6,j_2}N_{5,j_1}}{K_3}+2\begin{bmatrix} x_0 \\ u_0 \end{bmatrix}^T Q_{3,j} \\
     & g_{3,j}=\frac{d_{5,j_1}d_{6,j_2}-1}{K_3}-\begin{bmatrix} x_0 \\ u_0 \end{bmatrix}^T Q_{3,j}\begin{bmatrix} x_0 \\ u_0 \end{bmatrix} +H_{3,j}\begin{bmatrix} x_0 \\ u_0 \end{bmatrix} \\
    & N_{5,j_1}=K_3k_5n_{5,j_1}^T, ~d_{5,j_1}=K_3k_5b_{5,j_1}+1 \\
    & N_{6,j_2}=K_3k_6n_{6,j_2}^T, ~d_{6,j_2}=K_3k_6b_{6,j_2}+1 \\
    & C_{3,j}=\begin{bmatrix} R_{5,j_1} \\ R_{6,j_1} \end{bmatrix}  ~ c_{3,j}=\begin{bmatrix} r_{5,j_1} \\ r_{6,j_2} \end{bmatrix} + C_{3,j}\begin{bmatrix} x_0 \\ u_0 \end{bmatrix}.
   \end{split}
 \end{equation}
The numerical values of $Q_{i,j}, H_{i,j},g_{i,j},C_{i,j},c_{i,j}$, $i=2,3$, $j=1,\ldots,D_i$ can readily be computed from
the values $K_2,K_3,k_3,k_4,k_5,k_6$, $n_{k,j},b_{k,j}$, $k=3,4,5,6$, $j=1,\ldots,4$.  
and  \eqref{pl:approx:eq1}--\eqref{pl:approx:eq2}. 

We then apply Appendix \ref{lmi:appendix} to the functions $f_1,f_{2,pl}, f_{3,pl}$ as utility functions, i.e., the utility function of the first stakeholder will be $f_1$, and the utility function of the stakeholder $i$ will
be $f_{i,pl}$ for $i=1,2$. 
The resulting matrices $K$ and $Q$ are presented in Table \ref{num_result}. 
\begin{table}[H]
\caption{Strategy \eqref{safety_region} and safe set \eqref{safety_region} for the numerical example \label{num_result}}
\begin{tabular}{|c|c|}
\hline
    $K$ & \( \begin{bmatrix} 
0.0111  & 10.7842  & -6.7095 &   0.0393   & 1.5343
\end{bmatrix} \) \\
\hline
    $Q^{-1}$ & \(\begin{bmatrix} 
 0.0015  &  -0.1327   &-0.1364   &-0.0007  & -0.0439 \\
   -0.1327 & 246.5171 & -56.1857 &  -0.2685  &-17.5087 \\
   -0.1364  &-56.1857 & 247.0911 &  -0.2490 & -18.5745 \\ 
   -0.0007&   -0.2685   &-0.2490 &   0.0058  & -0.0812 \\
   -0.0439 & -17.5087  &-18.5745 &  -0.0812 &  26.8491
    \end{bmatrix} \) \\
 \hline
\end{tabular}
\end{table}

\subsection{Initial states and disturbances used for simulation}
\label{num_model:simul}
 In the dynamical model was simulated using different strategies, initial states, and was subjected to disturbances.
 Three different initial states were used: inside the safe set, outside the safe set, and around the equilibrium point. The numerical
 values of these initial states are described in Table \ref{initial_states:table}.
 \begin{table}[H]
 \caption{Initial states used for simulations \label{initial_states:table}}
 \begin{tabular}{|c|c|}
\hline
  Initial state $X(0)$ in $\mathcal{P}$ & $\begin{bmatrix} 99.23 &    0.16 &    0.16 &    43.93 &    12.99 \end{bmatrix}^T$ \\
\hline
  Initial state $X(0)$ is not in $\mathcal{P}$ & $\begin{bmatrix} 85 &    0.5 &    7.5 &   27.03 &   14.77 \end{bmatrix}^T$  \\
 \hline
Initial state $X(0)$ equals the equilibrium point $x_0$ & $\begin{bmatrix}  62.64 & 0.1 & 0.1 & 30.25 & 12.8    \end{bmatrix}^T$ \\
\hline  
\end{tabular}
\end{table}

 For the simulations, when  the disturbance $d(t)$ was not zero, it was chosen as follows:
  \begin{equation}
 \label{choice:disturbance1}
     d(t)=\left\{\begin{array}{rl}
                  -\sigma_x(Ax(t)+Bu(t))-\lambda d  & t=kN \\
                   0   & \mbox{ otherwise }
           \end{array}\right.
   \end{equation}
  i.e., $d(t)$ models a periodic change in the state , which occurs with a period $N$. 
  For the simulation we have chosen $\lambda=0.7$, $N=5$ and
  \begin{equation} 
 \label{choice:disturbance2}
      d=\begin{bmatrix} 40.3624   & 0.1000    &0.1000 &  20.3904  &  0.3000  &  0.1000 
\end{bmatrix}^T
   \end{equation}


\color{black}

\newpage



\end{document}